\documentclass[a4paper,fleqn]{cas-sc}

\usepackage{amsmath}
\usepackage{float}
\restylefloat{figure}
\usepackage{float}
\usepackage{placeins}
\usepackage{amsmath}
\usepackage{float}
\usepackage{subcaption}
\usepackage{tikz}
\usepackage{mathrsfs}
\usepackage{amsfonts}
\usepackage{amsmath}
\usepackage{array}
\usepackage{multirow}
\usepackage{longtable}
\usepackage{amsthm}
\usepackage{tabularx}
\usepackage{amssymb}
\usepackage{mathtools}
\usepackage{array,multirow}
\usepackage{xcolor}
\theoremstyle{definition}

\newtheorem{construction}{Construction}
\newtheorem{theorem}{Theorem}
\newtheorem{proposition}{Proposition}
\newtheorem{remark}{Remark}
\newtheorem{example}{Example}
\newtheorem{corollary}{Corollary}

\def\tsc#1{\csdef{#1}{\textsc{\lowercase{#1}}\xspace}}
\tsc{WGM}
\tsc{QE}

\begin{document}
	\let\WriteBookmarks\relax
	\def\floatpagepagefraction{1}
	\def\textpagefraction{.001}
	
	\shorttitle{}    
	
	\shortauthors{}  
	
	\title [mode = title]{Enumeration of Cospectral Graphs with Respect to the Average Flow Matrix}

	%

\author[1]{Usama Waheed}
\ead{usamawaheed@math.qau.edu.pk}
	
	\author[1,2]{Naveed Ahmed Azam\corref{cor1}}
	[
orcid=0000-0002-7941-3419]
\ead{azam@amp.i.kyoto-u.ac.jp}
\cortext[cor1]{Corresponding author}

\affiliation[1]{organization={Department of Mathematics, Quaid-i-Azam University},
city={Islamabad},
country={Pakistan}}

\affiliation[2]{organization={Graduate School of Advanced Integrated Studies in Human Survivability (Shishu-Kan), Kyoto University},
city={Kyoto},
postcode={606-8306},
country={Japan}}
	
	\begin{abstract}
		A graph is said to be determined by the spectrum of a graph-associated matrix if no non-isomorphic graph shares the same spectrum with it. In this paper, we introduce the \emph{average flow matrix} of a connected graph and study its spectral properties. We investigate its effectiveness in distinguishing graphs by carrying out a computational enumeration of connected graphs having a cospectral mate with respect to this matrix up to order $10$. Our computations show that every connected graph on at most $10$ vertices is determined by the spectrum of the average flow matrix. To the best of our knowledge, this is the first graph-associated matrix for which the spectrum distinguishes all connected graphs of order at most $10$.
	\end{abstract}
	

	
	\maketitle
	
	\section{Introduction}\label{sec1}
	Spectral graph theory studies the extent to which structural properties of a graph are reflected in the eigenvalues of graph-associated matrices. One of its central questions is whether a graph can be uniquely determined by the spectrum of such a matrix. In this context, graphs that share the same spectrum with a non-isomorphic graph are called cospectral, and the existence of such graphs shows that spectral data do not always determine graph structure completely. Nevertheless, spectral invariants provide useful tools for distinguishing graphs and are therefore also closely related to the graph isomorphism problem.

This question has a long history in spectral graph theory. Before 1956, it was believed that every graph is determined by its adjacency spectrum. G\"unthard and Primas~\cite{j} raised the broader question of whether there exists a graph-associated matrix for which every graph is determined by the corresponding spectrum. Later, \textcolor{black}{Godsil and McKay~\cite{r}} enumerated adjacency-cospectral graphs on up to $9$ vertices and introduced several construction techniques for cospectral graphs, among which Godsil--McKay switching is one of the best known. Brouwer and Spence~\cite{e} extended the enumeration of adjacency-cospectral graphs to graphs on $12$ vertices. Van Dam and Haemers~\cite{g} answered the question of G\"unthard and Primas~\cite{j} positively in principle by showing that it is possible to associate a matrix with a graph such that every graph is determined by its spectrum. However, no such concrete matrix has yet been identified~\cite{b}. Van Dam and Haemers~\cite{g} also conjectured that the proportion of graphs having a cospectral mate tends to $0$ as the order of the graph tends to infinity.

Motivated by these questions, a considerable amount of work has focused on the enumeration of cospectral graphs and on comparing the effectiveness of different graph-associated matrices. Haemers and Spence~\cite{f} studied enumeration with respect to the adjacency, Laplacian, and signless Laplacian matrices. Wilson and Zhu~\cite{y} showed that the fraction of graphs determined by spectral data can be improved by considering two or three matrices simultaneously. Aouchiche and Hansen~\cite{d} demonstrated that distance-based matrices perform better than the adjacency and Laplacian matrices in distinguishing graphs. More recently, Abiad and Alfaro~\cite{b} investigated both cospectral and coinvariant graphs and further refined the comparison of graph-associated matrices. These results suggest a common pattern: stronger enumeration results are often obtained when more structural information about a graph is incorporated into the associated matrix.

Despite these advances, none of the graph-associated matrices studied so far has been able to determine all connected graphs on $10$ vertices. Although distance-based matrices and combinations of several matrices lead to substantial improvements, they still leave non-isomorphic graphs that cannot be distinguished by the corresponding spectral or arithmetic invariants. This naturally motivates the search for new graph-associated matrices with stronger distinguishing power.

In this paper, we introduce the \emph{average flow matrix} of a graph. The diagonal entries of this matrix are derived from the average sums of the lengths of simple paths from a vertex to all other vertices, while the off-diagonal entries are determined by the distance structure of the graph. We enumerate cospectral graphs with respect to the average flow matrix and compare the resulting counts with those obtained for other graph-associated matrices. Our results show that every connected graph of order at most $10$ is determined by the spectrum of the average flow matrix. In particular, this provides a complete spectral characterization for connected graphs up to order $10$ and offers further evidence that almost all graphs may be determined by suitably chosen spectral data.

The paper is organized as follows. Section~\ref{sec2} contains the basic notation and terminology. Section~\ref{sec3} reviews the related work. In Section~\ref{sec4}, we introduce the average flow matrix, coflow graphs, and some theoretical results. Sections~\ref{sec5} and~\ref{sec6} present computational results. Finally, a conclusion is provided in Section~\ref{sec7}.

	\section{Preliminaries}\label{sec2}
\noindent Throughout this paper, we consider simple connected graphs. All paths are assumed to be simple unless stated otherwise.

Let $G=(V(G),E(G))$ be a graph, where $V(G)$ and $E(G)$ denote the vertex set and edge set of $G$, respectively. The \emph{order} and \emph{size} of $G$ are denoted by $|V(G)|$ and $|E(G)|$, respectively. For convenience, we write $n=|V(G)|$ and $m=|E(G)|$. For a vertex $v\in V(G)$, the \emph{degree} of $v$ is denoted by $d_G(v)$, and its \emph{neighbourhood} is defined by $N_G(v)=\{u\in V(G): uv\in E(G)\}$. A graph $G$ is called \emph{$r$-regular} if $d_G(v)=r$ for every $v\in V(G)$.

A \emph{simple path} in $G$ is a sequence of distinct vertices in which every two consecutive vertices are adjacent. The \emph{length} of a path is the number of edges it contains. For vertices $u,v\in V(G)$, the set of all simple paths joining $u$ and $v$ is denoted by $\mathcal{P}(u,v)$. The \emph{distance} between $u$ and $v$, denoted by $d_G(u,v)$, is the length of a shortest path joining them.

A graph $H$ is called a \emph{subgraph} of $G$ if $H\subseteq G$. Two graphs $G$ and $H$ are said to be \emph{isomorphic}, denoted by $G\cong H$, if there exists a bijection between their vertex sets that preserves adjacency. The \emph{complement} of $G$ is denoted by $\overline{G}$. A \emph{tree} is a connected graph containing no cycles.

The \emph{adjacency matrix} of a graph $G$ on vertices $v_1,v_2,\dots,v_n$ is denoted by $A(G)=(a_{ij})$, where $a_{ij}=1$ if $v_i$ and $v_j$ are adjacent, and $a_{ij}=0$ otherwise.

Let $\mathcal{G}_n$ denote the family of all connected graphs on $n$ vertices. For a real symmetric matrix $M$, the \emph{largest eigenvalue} and \emph{smallest eigenvalue} of $M$ are denoted by $\lambda_{\max}(M)$ and $\lambda_{\min}(M)$, respectively, and the \emph{trace} of $M$ is denoted by $\operatorname{tr}(M)$. The \emph{spread} of $M$ is defined by $\lambda_{\max}(M)-\lambda_{\min}(M)$. The diagonal matrix with diagonal entries $\ell_1,\ell_2,\dots,\ell_n$ is denoted by $\operatorname{diag}(\ell_1,\ell_2,\dots,\ell_n)$. We use $I_n$ and $J_n$ to denote the \emph{identity matrix} and the \emph{all-ones matrix} of order $n$, respectively. When the order is clear from the context, we simply write $I$ and $J$. The \emph{all-ones vector} of appropriate order is denoted by $\mathbf{1}$.

\section{Related Work}\label{sec3}

A graph is said to be \emph{determined by the spectrum} with respect to a given graph-associated matrix if no non-isomorphic graph shares the same spectrum with it. Spectral characterization problems of this kind have been studied for a variety of graph-associated matrices.

Harary and Read~\cite{t} gave one of the earliest enumerations of adjacency-cospectral graphs, determining the number of connected graphs on up to $7$ vertices that have an adjacency-cospectral mate.
The first major computational study in this direction is due to Godsil and McKay~\cite{r}, who enumerated adjacency-cospectral graphs on up to $9$ vertices. \textcolor{black}{In 2004}, Haemers and Spence~\cite{f} extended this line of work to the adjacency, Laplacian, and signless Laplacian matrices, and showed that the signless Laplacian matrix yields the best spectral discrimination among these three matrices. Let $G$ be a graph on vertices $v_1,v_2,\dots,v_n$, and let $\Delta(G)=\operatorname{diag}(d_G(v_1),d_G(v_2),\dots,d_G(v_n))$ be the \emph{degree matrix} of $G$. The \emph{Laplacian matrix} and \emph{signless Laplacian matrix} of $G$ are defined by $L(G)=\Delta(G)-A(G)$ and $Q(G)=\Delta(G)+A(G)$, respectively. Wilson and Zhu~\cite{y} later demonstrated that stronger discrimination can be obtained by combining several matrices simultaneously, namely the adjacency matrix $A$, the Laplacian matrix $L$, the signless Laplacian matrix $Q$, and the \emph{normalized Laplacian matrix} $\tilde{L}$, where $\tilde{L}(G)=I-\Delta(G)^{-1/2}A(G)\Delta(G)^{-1/2}$.

Distance-based matrices have proved to be particularly effective in the enumeration of cospectral graphs. Let $G$ be a connected graph on vertices $v_1,v_2,\dots,v_n$. The \emph{distance matrix} of $G$ is $D(G)=(d_G(v_i,v_j))_{i,j=1}^n$. The \emph{transmission} of a vertex $v_i$ is defined by $\operatorname{Tr}(v_i)=\sum_{j=1}^n d_G(v_i,v_j)$, and the corresponding diagonal matrix $\operatorname{Tr}(G)=\operatorname{diag}(\operatorname{Tr}(v_1),\operatorname{Tr}(v_2),\dots,\operatorname{Tr}(v_n))$ is called the \emph{transmission matrix} of $G$. The \emph{distance Laplacian matrix} and \emph{distance signless Laplacian matrix} of $G$ are defined by $D^L(G)=\operatorname{Tr}(G)-D(G)$ and $D^Q(G)=\operatorname{Tr}(G)+D(G)$, respectively. Aouchiche and Hansen~\cite{d} compared the adjacency, Laplacian, signless Laplacian, distance, distance Laplacian, and distance signless Laplacian matrices, and showed that the distance-based matrices provide substantially better discrimination. More recently, Abiad and Alfaro~\cite{b} studied both cospectral and coinvariant graphs, and obtained further improvements by combining adjacency- and distance-based invariants.

The present work continues this line of investigation by introducing the average flow matrix and studying its effectiveness in the enumeration of cospectral graphs.

\vspace{0.2cm}

We now recall some standard constructions that produce cospectral graphs.

\begin{construction}[\cite{c}, \emph{Godsil--McKay switching}]
	Let $G$ be a graph with vertex set $V(G)$, and let $\phi=(C_1,C_2,\dots,C_k,D)$ be a partition of $V(G)$ such that, for all $1 \leq i,j \leq k$ and for every vertex $v \in D$, the following conditions hold:
	\begin{itemize}
		\item[(i)] for any $v_i,v_j \in C_i$, $|N_{C_j}(v_i)|=|N_{C_j}(v_j)|$;
		\item[(ii)] the vertex $v$ is adjacent to either all, exactly half, or none of the vertices in $C_i$.
	\end{itemize}
	The graph $G'$ obtained from $G$ by Godsil--McKay switching is defined as follows: for each vertex $v \in D$, whenever $v$ is adjacent to exactly half of the vertices in some $C_i$, we delete those edges and join $v$ instead to the other half of the vertices in $C_i$.
\end{construction}

\begin{theorem}[\cite{c}]
	Let $G$ be a graph, and let $\phi$ be a partition of $V(G)$ satisfying the conditions of Godsil--McKay switching. Then the graph $G'$ obtained from $G$ by Godsil--McKay switching is cospectral with $G$, and their complements are also cospectral.
\end{theorem}

\begin{example}
	Consider a regular graph $G$, shown in Figure~\ref{cospectral_regular_graphs}, with an even number of vertices. Add a new vertex $v$ to $G$ and join it to exactly half of the vertices, namely $\{w_1,w_2,w_3,w_4,w_5\}$. Now construct a new graph $G'$ by deleting all edges incident with $v$ and then joining $v$ to the other half of the vertices, namely $\{v_1,v_2,v_3,v_4,v_5\}$. By Godsil--McKay switching, the graphs $G$ and $G'$ are cospectral.
	
	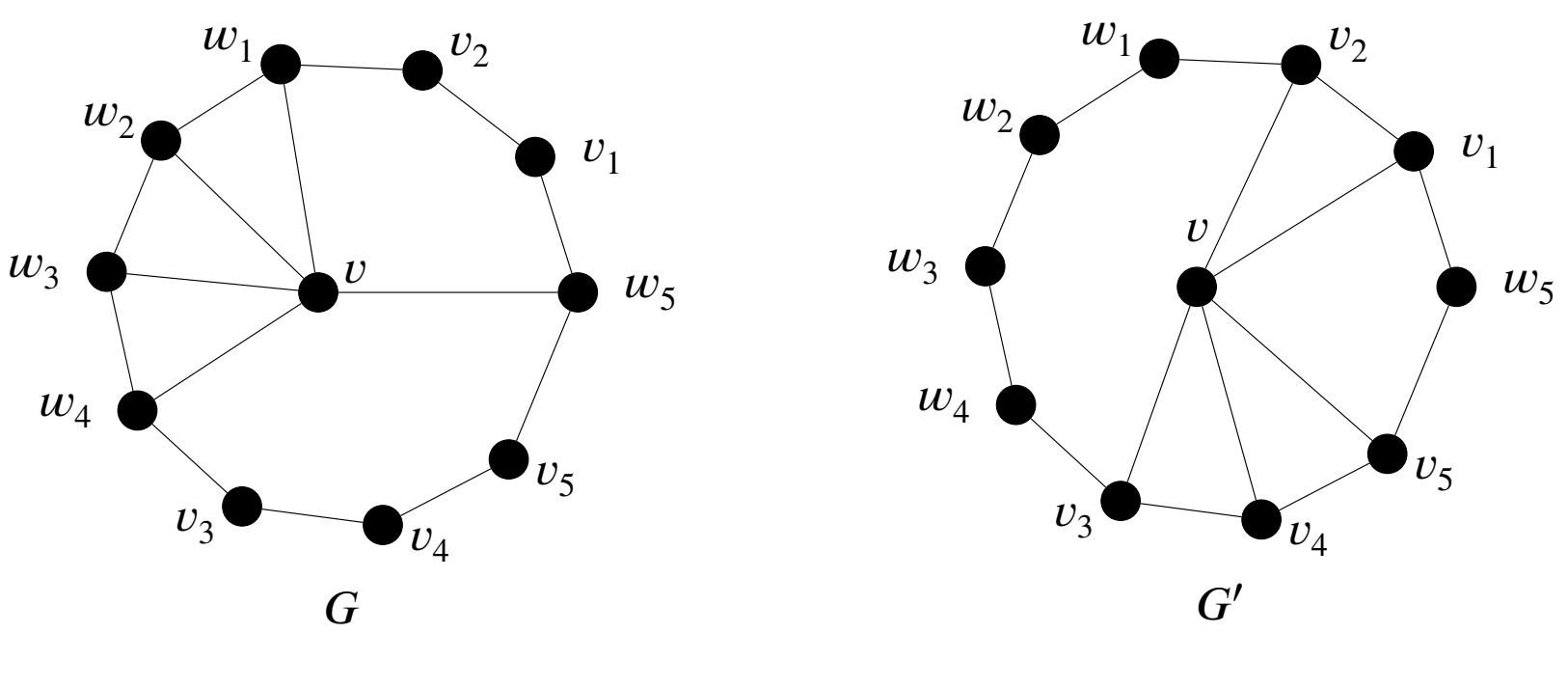
\begin{figure}[htbp]
		\centering
		\begin{minipage}[b]{0.3\textwidth}
			\centering
			\resizebox{\linewidth}{!}{%
				\begin{tikzpicture}[auto, every node/.style={circle, draw, fill=black, inner sep=5pt}]
					\node[fill=black, circle, label={[font=\LARGE, anchor=west]right:$v_1$}]  (0) at (35:3cm)  {};
					\node[fill=black, circle, label={[font=\LARGE, anchor=west]above:$v_2$}]  (1) at (70:3cm)  {};
					\node[fill=black, circle, label={[font=\LARGE, anchor=east]above:$w_1$}]  (2) at (105:3cm) {};
					\node[fill=black, circle, label={[font=\LARGE, anchor=east]above:$w_2$}]  (3) at (140:3cm) {};
					\node[fill=black, circle, label={[font=\LARGE, anchor=east]left:$w_3$}]   (4) at (175:3cm) {};
					\node[fill=black, circle, label={[font=\LARGE, anchor=east]left:$w_4$}]   (5) at (210:3cm) {};
					\node[fill=black, circle, label={[font=\LARGE, anchor=east]below:$v_3$}]  (6) at (245:3cm) {};
					\node[fill=black, circle, label={[font=\LARGE, anchor=west]below:$v_4$}]  (7) at (280:3cm) {};
					\node[fill=black, circle, label={[font=\LARGE, anchor=west]below:$v_5$}]  (8) at (315:3cm) {};
					\node[fill=black, circle, label={[font=\LARGE, anchor=west]right:$w_5$}]  (9) at (0:3cm)   {};
					\node[fill=black, circle, label={[font=\LARGE, anchor=west]above:$v$}]   (v) at (180:0.3cm) {};
					
					\draw (0)--(1)--(2)--(3)--(4)--(5)--(6)--(7)--(8)--(9)--(0);
					\draw (v) -- (2);
					\draw (v) -- (3);
					\draw (v) -- (4);
					\draw (v) -- (5);
					\draw (v) -- (9);
					\node[draw=none, fill=none] at (0,-4) {\LARGE$G$};
				\end{tikzpicture}%
			}
		\end{minipage}%
		\hspace{1cm}
		\begin{minipage}[b]{0.3\textwidth}
			\centering
			\resizebox{\linewidth}{!}{%
				\begin{tikzpicture}[auto, every node/.style={circle, draw, fill=black, inner sep=5pt}]
					\node[fill=black, circle, label={[font=\LARGE, anchor=west]right:$v_1$}]  (0) at (35:3cm)  {};
					\node[fill=black, circle, label={[font=\LARGE, anchor=west]above:$v_2$}]  (1) at (70:3cm)  {};
					\node[fill=black, circle, label={[font=\LARGE, anchor=east]above:$w_1$}]  (2) at (105:3cm) {};
					\node[fill=black, circle, label={[font=\LARGE, anchor=east]above:$w_2$}]  (3) at (140:3cm) {};
					\node[fill=black, circle, label={[font=\LARGE, anchor=east]left:$w_3$}]   (4) at (175:3cm) {};
					\node[fill=black, circle, label={[font=\LARGE, anchor=east]left:$w_4$}]   (5) at (210:3cm) {};
					\node[fill=black, circle, label={[font=\LARGE, anchor=east]below:$v_3$}]  (6) at (245:3cm) {};
					\node[fill=black, circle, label={[font=\LARGE, anchor=west]below:$v_4$}]  (7) at (280:3cm) {};
					\node[fill=black, circle, label={[font=\LARGE, anchor=west]below:$v_5$}]  (8) at (315:3cm) {};
					\node[fill=black, circle, label={[font=\LARGE, anchor=west]right:$w_5$}]  (9) at (0:3cm)   {};
					\node[fill=black, circle, label={[font=\LARGE, anchor=south]above:$v$}]   (v) at (180:0.3cm) {};
					
					\draw (0)--(1)--(2)--(3)--(4)--(5)--(6)--(7)--(8)--(9)--(0);
					\draw (v) -- (0);
					\draw (v) -- (1);
					\draw (v) -- (6);
					\draw (v) -- (7);
					\draw (v) -- (8);
					\node[draw=none, fill=none] at (0,-4) {\LARGE$G'$};
				\end{tikzpicture}%
			}
		\end{minipage}
		\caption{A pair of cospectral graphs obtained by Godsil--McKay switching.}
		\label{cospectral_regular_graphs}
	\end{figure}
\end{example}

\begin{construction}[\cite{c}, \emph{Construction via congruent binary matrices}]
	Let $A_1$ and $A_2$ be two matrices. We say that $A_1$ and $A_2$ are \emph{congruent} if $A_1^TA_1=A_2^TA_2$. Let $G$ be a graph on $n$ vertices, and let $K$ be an $m\times n$ binary matrix. We define $G(K)$ to be the graph with adjacency matrix
	\[
	\begin{bmatrix}
		0   & K    \\
		K^T & A(G)
	\end{bmatrix}.
	\]
\end{construction}

\begin{theorem}[\cite{c}]
	Let $H$ be a graph on $n$ vertices, and let $A$ and $B$ be binary congruent matrices of size $m\times n$. Then the graphs $H(A)$ and $H(B)$ are cospectral.
\end{theorem}

Let $A$ and $B$ be matrices of sizes $m\times n$ and $\ell\times k$, respectively. Their \emph{tensor product} $A \otimes B$ is the $(m\ell)\times (nk)$ matrix obtained by replacing each entry of $A$ by that entry multiplied by the entire matrix $B$.

\begin{construction}[\cite{c}, \emph{Tensor-product construction}]
	Let $G_1$ and $G_2$ be cospectral graphs, and let $X$ and $H$ be binary square matrices of the same order. Then the graph with adjacency matrix
	\[
	H \otimes I + X \otimes A(G_1)
	\]
	is cospectral with the graph having adjacency matrix
	\[
	H \otimes I + X \otimes A(G_2).
	\]
\end{construction}

\begin{construction}[\cite{c}, \emph{Tensor-product construction with complements}]
	Let $G_1$ and $G_2$ be cospectral graphs such that their complements $\bar{G}_1$ and $\bar{G}_2$ are also cospectral. Let $A_1, A_2, A_3,$ and $A_4$ be binary square matrices of the same order. Then the graph with adjacency matrix
	\[
	A_1\otimes I + A_2\otimes J + A_3\otimes A(G_1) + A_4\otimes A(\bar{G}_1)
	\]
	is cospectral with the graph having adjacency matrix
	\[
	A_1\otimes I + A_2\otimes J + A_3\otimes A(G_2) + A_4\otimes A(\bar{G}_2).
	\]
\end{construction}

\noindent Lepovi\'c~\cite{i} presented numerical data for graphs on $10$ vertices. \textcolor{black}{Extending the computation results of Godsil and McKay~\cite{r}}, Haemers and Spence~\cite{f} enumerated graphs on $10$ and $11$ vertices that are cospectral with respect to the adjacency, Laplacian, and signless Laplacian matrices. They also constructed such cospectral graphs by formulating Godsil--McKay switching as an operation on certain matrices.

\begin{theorem}[\cite{f}]
	Let $N$ be a binary matrix of size $b \times c$, and let $\tilde{N}$ be the matrix obtained from $N$ by Godsil--McKay switching. Let $B$ be a symmetric $b \times b$ matrix with constant row (and column) sums, and let $C$ be a symmetric $c \times c$ matrix. Define
	\[
	\begin{array}{cc}
		M =
		\begin{bmatrix}
			B & N \\
			N^T & C
		\end{bmatrix}
		&
		\tilde{M} =
		\begin{bmatrix}
			B & \tilde{N} \\
			\tilde{N}^T & C
		\end{bmatrix}.
	\end{array}
	\]
	Then $M$ and $\tilde{M}$ are cospectral.
\end{theorem}

The result of Haemers and Spence~\cite{f} gives a useful matrix formulation of Godsil--McKay switching and makes it possible to apply the method in a broader setting. In the same paper, Haemers and Spence~\cite{f} used this point of view in their computational study of cospectral graphs.

\begin{figure}
	\centering
	\begin{minipage}[b]{0.2\textwidth}
		\centering
		\resizebox{\linewidth}{!}{%
			\begin{tikzpicture}[auto, every node/.style={circle, draw, fill=black, inner sep=5pt}]
				\node[fill=black, circle, label={}] (0) at (90:3cm) {};
				\node[fill=black, circle, label={}] (1) at (180:3cm) {};
				\node[fill=black, circle, label={}] (2) at (270:3cm) {};
				\node[fill=black, circle, label={}] (3) at (0:3cm) {};
				
				\draw (0) -- (1) -- (2);
				\draw (0) -- (2);
			\end{tikzpicture}%
		}%
		\subcaption{}
	\end{minipage}
	\hspace{3cm}
	\begin{minipage}[b]{0.2\textwidth}
		\resizebox{\linewidth}{!}{%
			\begin{tikzpicture}[auto, every node/.style={circle, draw, fill=black, inner sep=5pt}]
				\node[fill=black, circle, label={}] (0) at (90:3cm) {};
				\node[fill=black, circle, label={}] (1) at (180:3cm) {};
				\node[fill=black, circle, label={}] (2) at (270:3cm) {};
				\node[fill=black, circle, label={}] (3) at (0:3cm) {};
				
				\draw (0) -- (1);
				\draw (0) -- (2);
				\draw (0) -- (3);
			\end{tikzpicture}%
		}%
		\subcaption{}
	\end{minipage}
	\caption{Laplacian-cospectral graphs with different numbers of triangles~\cite{f}.}
	\label{l_cospectral_graphs}
\end{figure}
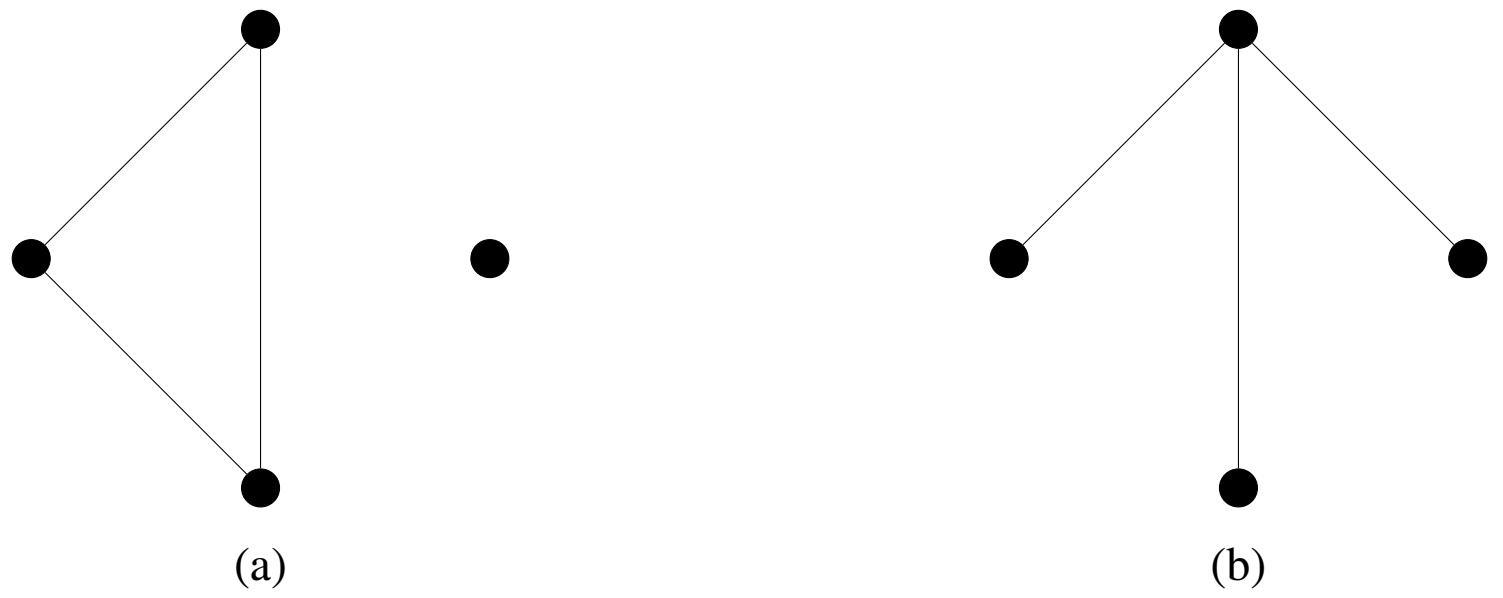

\noindent Since the number of graphs on $10$ vertices is $11{,}716{,}571$, the computation of cospectral graphs requires substantial time and memory. To deal with this, Haemers and Spence~\cite{f} exploited several necessary spectral conditions. For example, two graphs that are cospectral with respect to the adjacency matrix, and similarly with respect to the Laplacian matrix, must have the same number of edges. Moreover, adjacency-cospectral graphs must also have the same number of triangles. However, this is not true for Laplacian-cospectral graphs, as illustrated in Figure~\ref{l_cospectral_graphs}.

Brouwer and Spence~\cite{e} enumerated adjacency-cospectral graphs on $12$ vertices. They showed that the number of graphs having a cospectral mate initially increases and begins to decrease at $11$ vertices. \textcolor{black}{Wilson and Zhu~\cite{y}} enumerated cospectral graphs and trees with respect to the adjacency, Laplacian, signless Laplacian, and normalized Laplacian matrix $\tilde{L}$ of a graph $G$, as well as their combinations.

In 2016, Haemers~\cite{s} conjectured that almost all graphs are determined by their spectrum. He supported this conjecture with the following observations:
\begin{itemize}
	\item[(i)] The number of graphs having an adjacency-cospectral mate increases up to $10$ vertices, but begins to decrease at $11$ vertices~\cite{e}.
	
	\item[(ii)] Wang and Xu~\cite{u} considered the \emph{generalized spectrum}, that is, the spectrum of a graph $G$ together with the spectrum of its complement. Their experimental results indicate that a large proportion of the graphs tested are determined by their generalized spectrum.
	
	\item[(iii)] Wang~\cite{v,w} defined a class of graphs, denoted by $\mathcal{F}$, and proved that every graph in $\mathcal{F}$ is determined by its generalized spectrum. For a graph $G$ on $n$ vertices with adjacency matrix $A$, the \emph{walk matrix} $W$ is defined by
	\[
	W = [\mathbf{1}, A\mathbf{1}, A^{2}\mathbf{1}, \dots, A^{n-1}\mathbf{1}],
	\]
	where $\mathbf{1}$ denotes the all-ones vector. The class $\mathcal{F}$ consists of all graphs for which
	\[
	\frac{\det(W)}{2^{\lfloor n/2 \rfloor}}
	\]
	is an odd square-free integer, where $\det(W)$ denotes the \emph{determinant} of $W$.
\end{itemize}

Aouchiche and Hansen~\cite{d} enumerated graphs on $10$ vertices that are cospectral with respect to the distance matrix, the distance Laplacian matrix, and the distance signless Laplacian matrix. They also showed that these enumeration results can be sharpened by considering two or three of these matrices simultaneously, that is, by counting the number of graphs that have a cospectral mate with respect to each matrix in the chosen combination. Their best result for graphs on $10$ vertices was obtained by considering all three matrices together: only $7,622$ graphs were found to have a cospectral mate with respect to $(D, D^L, D^Q)$. Examples of graphs of smallest order having cospectral mates with respect to the associated distance matrices are shown in Figure~\ref{d_cosp_graphs}.

\begin{figure}[htbp]
	\centering
	\begin{minipage}{0.3\textwidth}
		\centering
		\resizebox{0.5\textwidth}{!}{%
			\begin{tikzpicture}[auto, every node/.style={circle, draw, fill=black, 
					inner sep=5pt}]
				\node[fill=black, circle, label={}] (0) at (51.3:3cm)  {};
				\node[fill=black, circle, label={}] (1) at (102.6:3cm) {};
				\node[fill=black, circle, label={}] (2) at (153.9:3cm) {};
				\node[fill=black, circle, label={}] (3) at (205.2:3cm) {};
				\node[fill=black, circle, label={}] (4) at (256.5:3cm) {};
				\node[fill=black, circle, label={}] (5) at (307.8:3cm) {};
				\node[fill=black, circle, label={}] (6) at (0:3cm)     {};
				\draw (0) -- (4); \draw (0) -- (5);
				\draw (1) -- (4); \draw (1) -- (5);
				\draw (2) -- (4); \draw (2) -- (6);
				\draw (3) -- (5); \draw (3) -- (6);
				\draw (4) -- (6); \draw (5) -- (6);
			\end{tikzpicture}%
		}
		\vspace{0.5cm}
		\resizebox{0.5\textwidth}{!}{%
			\begin{tikzpicture}[auto, every node/.style={circle, draw, fill=black, 
					inner sep=5pt}]
				\node[fill=black, circle, label={}] (0) at (51.3:3cm)  {};
				\node[fill=black, circle, label={}] (1) at (102.6:3cm) {};
				\node[fill=black, circle, label={}] (2) at (153.9:3cm) {};
				\node[fill=black, circle, label={}] (3) at (205.2:3cm) {};
				\node[fill=black, circle, label={}] (4) at (256.5:3cm) {};
				\node[fill=black, circle, label={}] (5) at (307.8:3cm) {};
				\node[fill=black, circle, label={}] (6) at (0:3cm)     {};
				\draw (0) -- (3); \draw (0) -- (5); \draw (0) -- (6);
				\draw (1) -- (4); \draw (1) -- (5); \draw (1) -- (6);
				\draw (2) -- (5); \draw (2) -- (6);
				\draw (3) -- (6); \draw (4) -- (6);
			\end{tikzpicture}%
		}
		\subcaption*{$D$-cospectral graphs.}
		\label{D-csp_graphs}
	\end{minipage}
	\begin{minipage}{0.3\textwidth}
		\centering
		\resizebox{0.5\textwidth}{!}{%
			\begin{tikzpicture}[auto, every node/.style={circle, draw, fill=black, 
					inner sep=5pt}]
				\node[fill=black, circle, label={}] (0) at (72:3cm)  {};
				\node[fill=black, circle, label={}] (1) at (144:3cm) {};
				\node[fill=black, circle, label={}] (2) at (216:3cm) {};
				\node[fill=black, circle, label={}] (3) at (288:3cm) {};
				\node[fill=black, circle, label={}] (4) at (0:3cm)   {};
				\draw (0) -- (3); \draw (0) -- (4);
				\draw (1) -- (3); \draw (1) -- (4);
				\draw (2) -- (3); \draw (2) -- (4);
				\draw (3) -- (4);
			\end{tikzpicture}%
		}
		\vspace{0.5cm}
		\resizebox{0.5\textwidth}{!}{%
			\begin{tikzpicture}[auto, every node/.style={circle, draw, fill=black, 
					inner sep=5pt}]
				\node[fill=black, circle, label={}] (0) at (72:3cm)  {};
				\node[fill=black, circle, label={}] (1) at (144:3cm) {};
				\node[fill=black, circle, label={}] (2) at (216:3cm) {};
				\node[fill=black, circle, label={}] (3) at (288:3cm) {};
				\node[fill=black, circle, label={}] (4) at (0:3cm)   {};
				\draw (0) -- (1); \draw (0) -- (2);
				\draw (0) -- (3); \draw (0) -- (4);
				\draw (2) -- (3); \draw (2) -- (4);
				\draw (3) -- (4);
			\end{tikzpicture}%
		}
		\subcaption*{$D^Q$-cospectral graphs.}
		\label{D^Q-csp_graphs}
	\end{minipage}
	\begin{minipage}{0.31\textwidth}
		\centering
		\resizebox{0.5\textwidth}{!}{%
			\begin{tikzpicture}[auto, every node/.style={circle, draw, fill=black, 
					inner sep=5pt}]
				\node[fill=black, circle, label={}] (0) at (51.3:3cm)  {};
				\node[fill=black, circle, label={}] (1) at (102.6:3cm) {};
				\node[fill=black, circle, label={}] (2) at (153.9:3cm) {};
				\node[fill=black, circle, label={}] (3) at (205.2:3cm) {};
				\node[fill=black, circle, label={}] (4) at (256.5:3cm) {};
				\node[fill=black, circle, label={}] (5) at (307.8:3cm) {};
				\node[fill=black, circle, label={}] (6) at (0:3cm)     {};
				\draw (0) -- (3); \draw (0) -- (2);
				\draw (0) -- (5); \draw (0) -- (6);
				\draw (1) -- (3); \draw (1) -- (4);
				\draw (1) -- (5); \draw (1) -- (6);
				\draw (2) -- (4); \draw (2) -- (5); \draw (2) -- (6);
				\draw (3) -- (5); \draw (3) -- (6);
				\draw (4) -- (6); \draw (5) -- (6);
			\end{tikzpicture}%
		}
		\vspace{0.5cm}
		\resizebox{0.5\textwidth}{!}{%
			\begin{tikzpicture}[auto, every node/.style={circle, draw, fill=black, 
					inner sep=5pt}]
				\node[fill=black, circle, label={}] (0) at (51.3:3cm)  {};
				\node[fill=black, circle, label={}] (1) at (102.6:3cm) {};
				\node[fill=black, circle, label={}] (2) at (153.9:3cm) {};
				\node[fill=black, circle, label={}] (3) at (205.2:3cm) {};
				\node[fill=black, circle, label={}] (4) at (256.5:3cm) {};
				\node[fill=black, circle, label={}] (5) at (307.8:3cm) {};
				\node[fill=black, circle, label={}] (6) at (0:3cm)     {};
				\draw (0) -- (3); \draw (0) -- (4);
				\draw (0) -- (5); \draw (0) -- (6);
				\draw (1) -- (3); \draw (1) -- (4); \draw (1) -- (5);
				\draw (2) -- (4); \draw (2) -- (5); \draw (2) -- (6);
				\draw (3) -- (4); \draw (3) -- (5); \draw (3) -- (6);
				\draw (4) -- (6); \draw (5) -- (6);
			\end{tikzpicture}%
		}
		\subcaption*{$D^L$-cospectral graphs.}
		\label{D^L-csp_graphs}
	\end{minipage}
	\centering
	\begin{minipage}{0.31\textwidth}
		\centering
		\resizebox{0.5\textwidth}{!}{%
			\begin{tikzpicture}[auto, every node/.style={circle, draw, fill=black, 
					inner sep=5pt}]
				\node[fill=black, circle, label={}] (0) at (45:3cm)  {};
				\node[fill=black, circle, label={}] (1) at (90:3cm)  {};
				\node[fill=black, circle, label={}] (2) at (135:3cm) {};
				\node[fill=black, circle, label={}] (3) at (180:3cm) {};
				\node[fill=black, circle, label={}] (4) at (225:3cm) {};
				\node[fill=black, circle, label={}] (5) at (270:3cm) {};
				\node[fill=black, circle, label={}] (6) at (315:3cm) {};
				\node[fill=black, circle, label={}] (7) at (360:3cm) {};
				\draw (0) -- (4); \draw (0) -- (6); \draw (0) -- (7);
				\draw (1) -- (5); \draw (1) -- (6); \draw (1) -- (7);
				\draw (2) -- (5); \draw (2) -- (6); \draw (2) -- (7);
				\draw (3) -- (5); \draw (3) -- (7);
				\draw (4) -- (7); \draw (5) -- (7); \draw (6) -- (7);
			\end{tikzpicture}%
		}
		\vspace{1cm}
		\resizebox{0.5\textwidth}{!}{%
			\begin{tikzpicture}[auto, every node/.style={circle, draw, fill=black, 
					inner sep=5pt}]
				\node[fill=black, circle, label={}] (0) at (45:3cm)  {};
				\node[fill=black, circle, label={}] (1) at (90:3cm)  {};
				\node[fill=black, circle, label={}] (2) at (135:3cm) {};
				\node[fill=black, circle, label={}] (3) at (180:3cm) {};
				\node[fill=black, circle, label={}] (4) at (225:3cm) {};
				\node[fill=black, circle, label={}] (5) at (270:3cm) {};
				\node[fill=black, circle, label={}] (6) at (315:3cm) {};
				\node[fill=black, circle, label={}] (7) at (360:3cm) {};
				\draw (0) -- (4); \draw (0) -- (7);
				\draw (1) -- (4); \draw (1) -- (7);
				\draw (2) -- (5); \draw (2) -- (6); \draw (2) -- (7);
				\draw (3) -- (5); \draw (3) -- (6); \draw (3) -- (7);
				\draw (4) -- (6); \draw (4) -- (7);
				\draw (5) -- (7); \draw (6) -- (7);
			\end{tikzpicture}%
		}
		\subcaption*{$(D^L, D^Q)$-cospectral graphs.}
	\end{minipage}
	\begin{minipage}{0.31\textwidth}
		\centering
		\resizebox{0.5\textwidth}{!}{%
			\begin{tikzpicture}[auto, every node/.style={circle, draw, fill=black, 
					inner sep=5pt}]
				\node[fill=black, circle, label={}] (0) at (36:3cm)  {};
				\node[fill=black, circle, label={}] (1) at (72:3cm)  {};
				\node[fill=black, circle, label={}] (2) at (108:3cm) {};
				\node[fill=black, circle, label={}] (3) at (144:3cm) {};
				\node[fill=black, circle, label={}] (4) at (180:3cm) {};
				\node[fill=black, circle, label={}] (5) at (216:3cm) {};
				\node[fill=black, circle, label={}] (6) at (252:3cm) {};
				\node[fill=black, circle, label={}] (7) at (288:3cm) {};
				\node[fill=black, circle, label={}] (8) at (324:3cm) {};
				\node[fill=black, circle, label={}] (9) at (0:3cm)   {};
				\draw (0) -- (3); \draw (0) -- (5);
				\draw (0) -- (7); \draw (0) -- (9);
				\draw (1) -- (4); \draw (1) -- (6);
				\draw (1) -- (7); \draw (1) -- (8);
				\draw (2) -- (3); \draw (2) -- (5);
				\draw (2) -- (6); \draw (2) -- (7); \draw (2) -- (8);
				\draw (3) -- (6); \draw (3) -- (9);
				\draw (4) -- (7); \draw (4) -- (8); \draw (4) -- (9);
				\draw (5) -- (7); \draw (5) -- (8);
				\draw (6) -- (9);
			\end{tikzpicture}%
		}
		\vspace{1cm}
		\resizebox{0.5\textwidth}{!}{%
			\begin{tikzpicture}[auto, every node/.style={circle, draw, fill=black, 
					inner sep=5pt}]
				\node[fill=black, circle, label={}] (0) at (36:3cm)  {};
				\node[fill=black, circle, label={}] (1) at (72:3cm)  {};
				\node[fill=black, circle, label={}] (2) at (108:3cm) {};
				\node[fill=black, circle, label={}] (3) at (144:3cm) {};
				\node[fill=black, circle, label={}] (4) at (180:3cm) {};
				\node[fill=black, circle, label={}] (5) at (216:3cm) {};
				\node[fill=black, circle, label={}] (6) at (252:3cm) {};
				\node[fill=black, circle, label={}] (7) at (288:3cm) {};
				\node[fill=black, circle, label={}] (8) at (324:3cm) {};
				\node[fill=black, circle, label={}] (9) at (0:3cm)   {};
				\draw (0) -- (3); \draw (0) -- (5);
				\draw (0) -- (7); \draw (0) -- (8);
				\draw (1) -- (4); \draw (1) -- (6);
				\draw (1) -- (7); \draw (1) -- (8);
				\draw (2) -- (4); \draw (2) -- (5);
				\draw (2) -- (6); \draw (2) -- (7); \draw (2) -- (8);
				\draw (3) -- (6); \draw (3) -- (8); \draw (3) -- (9);
				\draw (4) -- (7); \draw (4) -- (8);
				\draw (5) -- (8); \draw (5) -- (9);
				\draw (6) -- (9);
			\end{tikzpicture}%
		}
		\subcaption*{$(D, D^L, D^Q)$-cospectral graphs.}
	\end{minipage}
	\caption{Cospectral graphs with respect to distance matrices and their combinations.}
	\label{d_cosp_graphs}
\end{figure}

Abiad and Alfaro~\cite{b} enumerated cospectral and coinvariant graphs, that is, graphs having the same \emph{Smith normal form}, for orders up to $10$. Let $M$ be an integer matrix of rank $r$. The \emph{Smith normal form} of $M$ is the diagonal matrix
\[
\operatorname{diag}(d_1,d_2,\dots,d_r,0,\dots,0),
\]
where each $d_i$ is a positive integer and $d_i$ divides $d_{i+1}$ for $1 \leq i < r$, obtained from $M$ by multiplication on the left and right by unimodular integer matrices. Two graphs are said to be \emph{coinvariant} with respect to a graph-associated matrix if the corresponding matrices have the same Smith normal form. They showed that the Smith normal form is not a particularly effective invariant for the adjacency, Laplacian, and signless Laplacian matrices, but performs much better for distance-based matrices~\cite{b}. For example, on $10$ vertices there are $787{,}851$ graphs having a $D^L$-cospectral mate, whereas only $642{,}002$ graphs have a $D^L$-coinvariant mate. Similarly, there are $319{,}324$ graphs having a $D^Q$-cospectral mate, whereas only $264{,}955$ graphs have a $D^Q$-coinvariant mate.

Following the approach in~\cite{y,d}, they further refined their enumeration by considering combinations of matrices. In~\cite{d}, it was shown that there are $7{,}622$ graphs on $10$ vertices having a cospectral mate with respect to $(D, D^L, D^Q)$. In~\cite{b}, it was shown that there are $7{,}490$ graphs on $10$ vertices having an adjacency-cospectral mate that is also a coinvariant mate with respect to $D^L$ and $D^Q$.

For ease of comparison, we use the following notation in Table~\ref{table:important}. For a graph-associated matrix $M$, let $\mathcal{G}_n^{cs}(M)$ denote the family of connected graphs on $n$ vertices that have a cospectral mate with respect to $M$, and let $\mathcal{G}_n^{in}(M)$ denote the family of connected graphs on $n$ vertices that have a coinvariant mate with respect to $M$. More generally, for graph-associated matrices $M_1,M_2,\dots,M_k$, let $\mathcal{G}_n^{cs}(M_1,M_2,\dots,M_k)$ denote the family of connected graphs on $n$ vertices having a common cospectral mate with respect to each of the matrices $M_1,M_2,\dots,M_k$. Mixed notation is interpreted componentwise; for example, $\mathcal{G}_n^{cs\text{-}in}(M_1,M_2)$ denotes the family of connected graphs on $n$ vertices having a mate that is cospectral with respect to $M_1$ and coinvariant with respect to $M_2$. The enumeration results from the articles discussed above are summarized in Table~\ref{table:important}.

\begin{table}[ht]
	\caption{Number of graphs having a cospectral or coinvariant mate for various families of connected graphs.\label{table:important}}
	\begin{tabular*}{\tblwidth}{@{}p{4cm}p{4.4cm}cccccc@{}}
		\toprule
		& & \multicolumn{6}{c}{$n$} \\
		\cmidrule(lr){3-8}
		\textbf{References} & \textbf{Family} & \textbf{5} & \textbf{6} & \textbf{7} & \textbf{8} & \textbf{9} & \textbf{10} \\
		\midrule
		& $\mathcal{G}_{n}$ 
		& 21 & 112 & 853 & 11,117 & 261,080 & 11,716,571 \\
		\cite{t} 
		& $\mathcal{G}_{n}^{cs}(A) $ 
		& 0 & 2 & 63 & NA & NA & NA \\
		\cite{r} 
		& $\mathcal{G}_{n}^{cs}(A)$ 
		& 0 & 2 & 63 & 1,353 & 46,930 & NA \\
		\multirow{3}{*}{\cite{f}}
		& $\mathcal{G}_{n}^{cs}(A)$ 
		& 0 & 2 & 63 & 1,353 & 46,930 & 2,462,141 \\
		& $\mathcal{G}_{n}^{cs}(L)$ 
		& 0 & 4 & 115 & 1,611 & 40,560 & 1,367,215 \\
		& $\mathcal{G}_{n}^{cs}(Q)$ 
		& 2 & 10 & 80 & 1,047 & 17,627 & 615,919 \\[4pt]
		\multirow{6}{*}{\cite{y}}
		& $\mathcal{G}_{n}^{cs}(A, L)$ 
		& 0 & 0 & 0 & 0 & 82 & 13,864 \\
		& $\mathcal{G}_{n}^{cs}(L, Q)$ 
		& 0 & 0 & 16 & 232 & 4,139 & 107,835 \\
		& $\mathcal{G}_{n}^{cs}(A, Q)$ 
		& 0 & 0 & 0 & 0 & 8 & 10,716 \\
		& $\mathcal{G}_{n}^{cs}(A, \tilde{L})$ 
		& 0 & 0 & 0 & 6 & 14 & 10,281 \\
		& $\mathcal{G}_{n}^{cs}(L, \tilde{L})$ 
		& 0 & 0 & 0 & 0 & 6 & 10,256 \\
		& $\mathcal{G}_{n}^{cs}(A, L, Q, \tilde{L})$ 
		& 0 & 0 & 0 & 0 & 2 & 10,124 \\[4pt]
		\multirow{7}{*}{\cite{d}}
		& $\mathcal{G}_{n}^{cs}(D)$ 
		& 0 & 0 & 22 & 658 & 25,058 & 1,389,986 \\
		& $\mathcal{G}_{n}^{cs}(D^L)$ 
		& 0 & 0 & 43 & 745 & 20,455 & 787,851 \\
		& $\mathcal{G}_{n}^{cs}(D^Q)$ 
		& 2 & 6 & 38 & 453 & 8,168 & 319,324 \\
		& $\mathcal{G}_{n}^{cs}(D, D^L)$ 
		& 0 & 0 & 0 & 0 & 32 & 9,449 \\
		& $\mathcal{G}_{n}^{cs}(D, D^Q)$ 
		& 0 & 0 & 0 & 0 & 0 & 7,712 \\
		& $\mathcal{G}_{n}^{cs}(D^L, D^Q)$ 
		& 0 & 0 & 0 & 90 & 1,965 & 61,909 \\
		& $\mathcal{G}_{n}^{cs}(D, D^L, D^Q)$ 
		& 0 & 0 & 0 & 0 & 0 & 7,622 \\[4pt]
		\multirow{18}{*}{\cite{b}}
		& $\mathcal{G}_{n}^{in}(A)$ 
		& 20 & 112 & 853 & 11,117 & 261,080 & 11,716,571 \\
		& $\mathcal{G}_{n}^{in}(L)$ 
		& 8 & 57 & 526 & 8,027 & 221,834 & 11,036,261 \\
		& $\mathcal{G}_{n}^{in}(Q)$ 
		& 11 & 78 & 620 & 7,962 & 201,282 & 10,086,812 \\
		& $\mathcal{G}_{n}^{in}(D)$ 
		& 15 & 102 & 835 & 11,080 & 260,991 & 11,716,249 \\
		& $\mathcal{G}_{n}^{in}(D^L)$ 
		& 0 & 0 & 18 & 455 & 16,505 & 642,002 \\
		& $\mathcal{G}_{n}^{in}(D^Q)$ 
		& 2 & 4 & 20 & 259 & 7,444 & 264,955 \\
		& $\mathcal{G}_{n}^{in}(D, D^L)$ 
		& 0 & 0 & 0 & 32 & 1,770 & 92,915 \\
		& $\mathcal{G}_{n}^{in}(D, D^Q)$ 
		& 0 & 0 & 0 & 20 & 432 & 24,517 \\
		& $\mathcal{G}_{n}^{in}(D^L, D^Q)$ 
		& 0 & 0 & 0 & 44 & 1,447 & 46,239 \\
		& $\mathcal{G}_{n}^{cs-in}(A, D^L)$ 
		& 0 & 0 & 0 & 0 & 32 & 9,253 \\
		& $\mathcal{G}_{n}^{cs-in}(D, D^L)$ 
		& 0 & 0 & 0 & 0 & 32 & 8,935 \\
		& $\mathcal{G}_{n}^{cs-in}(A, D^Q)$ 
		& 0 & 0 & 0 & 0 & 2 & 7,688 \\
		& $\mathcal{G}_{n}^{cs-in}(D, D^Q)$ 
		& 0 & 0 & 0 & 0 & 0 & 7,646 \\
		& $\mathcal{G}_{n}^{in}(D, D^L, D^Q)$ 
		& 0 & 0 & 0 & 0 & 138 & 12,246 \\
		& $\mathcal{G}_{n}^{cs-cs-in}(A, D, D^L)$ 
		& 0 & 0 & 0 & 0 & 32 & 8,743 \\
		& $\mathcal{G}_{n}^{cs-cs-in}(A, D, D^Q)$ 
		& 0 & 0 & 0 & 0 & 0 & 7,550 \\
		& $\mathcal{G}_{n}^{cs-in-in}(A, D^Q, D^L)$ 
		& 0 & 0 & 0 & 0 & 0 & 7,490 \\
		& $\mathcal{G}_{n}^{cs-in-in}(D, D^L, D^Q)$ 
		& 0 & 0 & 0 & 0 & 0 & 7,510 \\
		\bottomrule
	\end{tabular*}
\end{table}

Table~\ref{table:important} shows that the use of several matrices together substantially reduces the number of graphs having a cospectral or coinvariant mate. In particular, the strongest reductions are obtained when distance-based matrices are combined with additional spectral or arithmetic information. Nevertheless, even these refinements still leave a nontrivial number of unresolved cases. This naturally leads to the search for other graph-associated matrices with stronger distinguishing power, which motivates the study undertaken in the present paper.

\newpage	
	\section{The Average Flow Matrix of a Graph}\label{sec4}
	
We study graph-associated matrices and their spectra in order to improve the enumeration of
cospectral graphs. Existing results indicate that stronger spectral discrimination is often
achieved by incorporating richer structural information into the associated
matrix~\cite{b,d,f,r}. For instance, the Laplacian matrix~\cite{f} distinguishes more
graphs than the adjacency matrix~\cite{r}, since it reflects both adjacency and degree
information. Similarly, the distance matrix~\cite{d} outperforms both the adjacency and the
Laplacian matrices~\cite{d,f}, as it encodes the pairwise distances between vertices.
Furthermore, the distance Laplacian and distance signless Laplacian matrices~\cite{d}
provide an additional improvement over the distance matrix. Motivated by these
observations, we introduce the \emph{average flow matrix} of a graph, which augments the
negative distance matrix with a vertex-wise quantity measuring the average length of all
simple paths originating at each vertex.

Let $G$ be a connected graph with vertex set $V(G)=\{v_1,\dots,v_n\}$, and let $v\in V(G)$.
For each vertex $u\neq v$, let $\mathcal{P}(v,u)$ denote the set of all simple paths from
$v$ to $u$. The \emph{average flow} of $v$, denoted by $\mathrm{fl}(v)$, is defined by
\begin{equation}\label{eq:flow}
	\mathrm{fl}(v) \;\triangleq\;
	\frac{\displaystyle\sum_{u\neq v}\ \sum_{P\in\mathcal{P}(v,u)} |P|}{p_v},
	\qquad
	p_v \;=\; \sum_{u\neq v}\bigl|\mathcal{P}(v,u)\bigr|,
\end{equation}
where $|P|$ denotes the length of the path $P$. Thus, $\mathrm{fl}(v)$ is the total length
of all simple paths from $v$ to all other vertices, divided by the total number $p_v$ of
such paths. Since $G$ is finite, each set $\mathcal{P}(v,u)$ is finite and nonempty, and
therefore $\mathrm{fl}(v)$ is a well-defined positive rational number. Moreover,
$\mathrm{fl}(v)$ is invariant under graph isomorphism.

Using these vertex average flows, we define the associated matrix as follows. The
\emph{average flow matrix} $F(G)=[a_{ij}]$ of $G$ is the $n\times n$ symmetric matrix given
by
\begin{equation}\label{eq:flowmatrix}
	a_{ii}=\mathrm{fl}(v_i),\qquad
	a_{ij}=-\,\operatorname{dist}_G(v_i,v_j)\quad(i\neq j),
\end{equation}
where $\operatorname{dist}_G(v_i,v_j)$ denotes the length of a shortest path between
$v_i$ and $v_j$ in $G$. Equivalently,
\[
	F(G)=\operatorname{diag}\bigl(\mathrm{fl}(v_1),\dots,\mathrm{fl}(v_n)\bigr)-D(G),
\]
where $D(G)$ is the distance matrix of $G$.


Figure~\ref{graph_with_flow_matrix} illustrates a graph $G$ together with its average flow
matrix $F(G)$.

\begin{figure}[htbp]
	\centering
	\begin{minipage}[c]{0.3\textwidth}
		\centering
		\resizebox{0.8\linewidth}{!}{%
			\begin{tikzpicture}[auto, every node/.style={circle, draw, fill=black,
					inner sep=5pt}]
				\node[fill=black, circle,
				label={[font=\LARGE]above:$v_2$}] (0) at (60:3cm)  {};
				\node[fill=black, circle,
				label={[font=\LARGE]above:$v_1$}] (1) at (120:3cm) {};
				\node[fill=black, circle,
				label={[font=\LARGE]left:$v_6$}]  (2) at (180:3cm) {};
				\node[fill=black, circle,
				label={[font=\LARGE]below:$v_5$}] (3) at (240:3cm) {};
				\node[fill=black, circle,
				label={[font=\LARGE]below:$v_4$}] (4) at (300:3cm) {};
				\node[fill=black, circle,
				label={[font=\LARGE]right:$v_3$}] (5) at (0:3cm)   {};

				\draw (0) -- (3);
				\draw (0) -- (4);
				\draw (0) -- (5);
				\draw (1) -- (3);
				\draw (1) -- (5);
				\draw (2) -- (4);
				\draw (2) -- (5);
				\node[draw=none, fill=none] at (0,-4) {\Large$G$};
			\end{tikzpicture}%
		}
	\end{minipage}%
	\hspace{0.001\textwidth}%
	\begin{minipage}[c]{0.60\textwidth}
		\centering
		$F(G)=
		\begin{bmatrix}
			2.7333 & -2     & -2     & -1     & -1     & -1     \\
			-2     & 3.1764 & -2     & -1     & -3     & -1     \\
			-2     & -2     & 3.1764 & -3     & -1     & -1     \\
			-1     & -1     & -3     & 3.1764 & -2     & -2     \\
			-1     & -3     & -1     & -2     & 3.1764 & -2     \\
			-1     & -1     & -1     & -2     & -2     & 2.7333
		\end{bmatrix}$
	\end{minipage}
	\caption{A graph $G$ together with its average flow matrix $F(G)$. The diagonal entries
	are rounded to four decimal places.}
	\label{graph_with_flow_matrix}
\end{figure}
We next establish a simple relationship between the distance spectrum of a graph and the spectrum of its average flow matrix. We say that a graph $G$ is \emph{flow-stable} if all its vertices have the same average flow, that is,
\[
	\mathrm{fl}(v_1)=\mathrm{fl}(v_2)=\cdots=\mathrm{fl}(v_n).
\]

\begin{proposition}\label{prop:flow-stable}
	Let $G$ be a flow-stable graph of order $n$, let $c$ denote the common average flow of
	its vertices, and let $D(G)$ be its distance matrix with distance spectrum
	$\{d_1,d_2,\dots,d_n\}$, listed with multiplicities. Then the average flow matrix $F(G)$
	has spectrum
	\[
		\{\,c-d_1,\ c-d_2,\ \dots,\ c-d_n\,\},
	\]
	with the same multiplicities. Moreover, $F(G)$ and $D(G)$ have the same eigenvectors.
\end{proposition}

\begin{proof}
	Since $G$ is flow-stable, every diagonal entry of $F(G)$ equals $c$. Hence, by
	\eqref{eq:flowmatrix},
	\[
		F(G)=cI-D(G).
	\]
	If $x$ is an eigenvector of $D(G)$ corresponding to the eigenvalue $d_i$, then
	\[
		F(G)x=(cI-D(G))x=cx-d_ix=(c-d_i)x.
	\]
	Therefore, $x$ is an eigenvector of $F(G)$ corresponding to the eigenvalue $c-d_i$.
	It follows that the spectrum of $F(G)$ is $\{c-d_i\}_{i=1}^n$, with the same
	multiplicities, and that $F(G)$ and $D(G)$ have the same eigenvectors.
\end{proof}
	As a direct application of Proposition~\ref{prop:flow-stable}, we obtain the following
corollary for vertex-transitive graphs.

\begin{corollary}\label{cor:vertex-transitive}
	Every vertex-transitive graph is flow-stable. Consequently, if $G$ is a
	vertex-transitive graph with common average flow $c$ and distance spectrum
	$\{d_1,\dots,d_n\}$, then the spectrum of its average flow matrix is
	$\{c-d_i\}_{i=1}^{n}$.
\end{corollary}

\begin{proof}
	If $G$ is vertex-transitive, then for any two vertices $u,w\in V(G)$ there exists an
	automorphism of $G$ mapping $u$ to $w$. Since $\mathrm{fl}(v)$ is invariant under graph
	automorphisms, it follows that all vertices of $G$ have the same average flow. Hence $G$
	is flow-stable, and the conclusion follows from Proposition~\ref{prop:flow-stable}.
\end{proof}

%
We now determine the spectrum of the average flow matrix of the complete graph $K_n$.

\begin{proposition}\label{prop:Kn-spectrum}
	Let $K_n$ be the complete graph on $n\geq 2$ vertices, and let $c=\mathrm{fl}(v)$ denote
	the common average flow of its vertices. Then
	\[
		\operatorname{spec}(F(K_n))
		=
		\{\, (c+1)^{[n-1]},\ c-(n-1)\,\}.
	\]
	In particular, the largest eigenvalue of $F(K_n)$ is
	\[
		\lambda_{\max}(F(K_n))=c+1.
	\]
\end{proposition}

\begin{proof}
	Since $K_n$ is vertex-transitive, it is flow-stable, and hence every diagonal entry of
	$F(K_n)$ equals $c$. Moreover, every pair of distinct vertices in $K_n$ is adjacent, so
	$\operatorname{dist}_{K_n}(u,v)=1$ whenever $u\neq v$. Therefore,
	\[
		F(K_n)=cI-D(K_n).
	\]
	It is well known that
	\[
		\operatorname{spec}(D(K_n))=\{\,n-1,\ (-1)^{[n-1]}\,\}.
	\]
	By Proposition~\ref{prop:flow-stable}, it follows that
	\[
		\operatorname{spec}(F(K_n))
		=
		\{\,c-(n-1),\ (c+1)^{[n-1]}\,\}.
	\]
	Since $n\geq 2$, we have $c-(n-1)<c+1$, and therefore
	\[
		\lambda_{\max}(F(K_n))=c+1.
	\]
\end{proof}
\begin{corollary}\label{cor:Kn-lambda-max}
	Let $\lambda_{\max}(F(K_n))$ denote the largest eigenvalue of the average flow matrix of
	$K_n$. If $n-2\leq \mathrm{fl}(v)\leq n-1$, then
	\[
		n-1\leq \lambda_{\max}(F(K_n))\leq n.
	\]
\end{corollary}

\begin{proof}
	By Proposition~\ref{prop:Kn-spectrum},
	\[
		\lambda_{\max}(F(K_n))=\mathrm{fl}(v)+1.
	\]
	The conclusion follows immediately from the bound
	$n-2\leq \mathrm{fl}(v)\leq n-1$.
\end{proof}

\subsection{Coflow Graphs}\label{subsec:coflow}

Enumerating cospectral graphs of larger order is computationally demanding because the
number of connected graphs grows rapidly with the number of vertices. For example, there are
$11{,}716{,}571$ connected graphs on 10 vertices. A direct pairwise spectral comparison over
such a large collection is prohibitively expensive on standard hardware. To reduce the
search space, we first identify a larger class of graphs that necessarily contains all
cospectral pairs.

For a connected graph $G$, we define its \emph{average flow sum} by
\begin{equation}\label{eq:flow-sum}
	\mathrm{fl}(G)\triangleq \sum_{v_i\in V(G)} \mathrm{fl}(v_i).
\end{equation}
Two graphs are said to be \emph{coflow mates} if they have the same average flow sum.
Equivalently, graphs with the same value of $\mathrm{fl}(G)$ belong to the same
\emph{coflow class}.

Since the diagonal entries of the average flow matrix $F(G)$ are precisely the vertex
average flows, we have
\[
	\operatorname{tr}(F(G))=\mathrm{fl}(G).
\]
Moreover, the trace of a square matrix is equal to the sum of its eigenvalues. It follows
that if two graphs are cospectral with respect to the average flow matrix, then they must
have the same average flow sum. Therefore,
\[
	\mathrm{fl}(G)\neq \mathrm{fl}(H)
	\qquad\Longrightarrow\qquad
	G \text{ and } H \text{ are not cospectral with respect to } F.
\]
Thus, coflow classes provide a natural and efficient first filter in the enumeration of
cospectral graphs.

We enumerated the connected graphs belonging to non-singleton coflow classes for all orders
up to 10. The results are presented in Table~\ref{tab2}. Only a small fraction of connected
graphs have a coflow mate; for example, among the $11{,}716{,}571$ connected graphs on 10
vertices, only $730{,}276$ belong to non-singleton coflow classes, which is approximately
$ 6.23\%$.

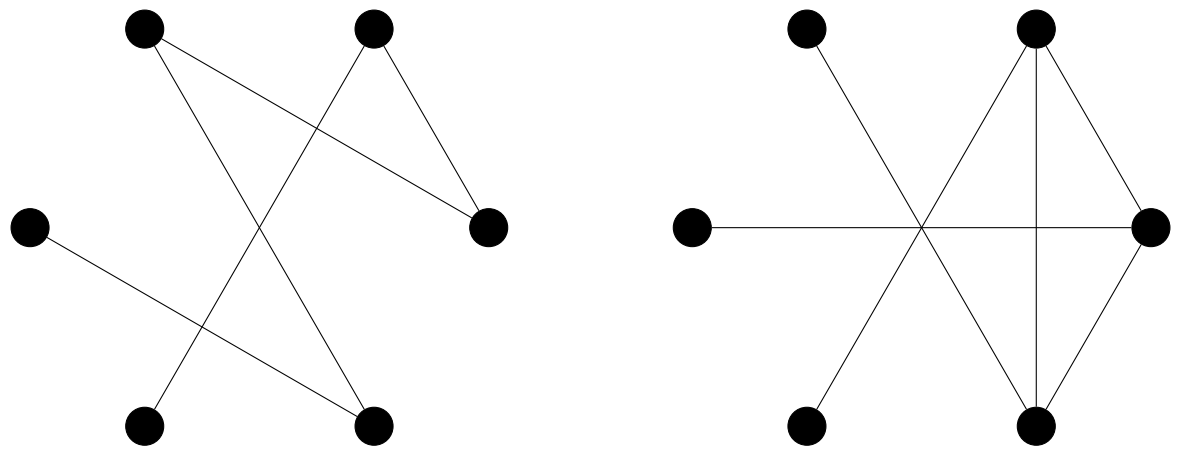
\begin{figure}[htbp]
	\centering
	\begin{minipage}[b]{0.2\textwidth}
		\centering
		\resizebox{\linewidth}{!}{%
			\begin{tikzpicture}[auto, every node/.style={circle, draw, fill=black,
					inner sep=5pt}]
				\node[fill=black, circle, label={}] (0) at (60:3cm)  {};
				\node[fill=black, circle, label={}] (1) at (120:3cm) {};
				\node[fill=black, circle, label={}] (2) at (180:3cm) {};
				\node[fill=black, circle, label={}] (3) at (240:3cm) {};
				\node[fill=black, circle, label={}] (4) at (300:3cm) {};
				\node[fill=black, circle, label={}] (5) at (0:3cm)   {};

				\draw (0) -- (3);
				\draw (0) -- (5);
				\draw (1) -- (4);
				\draw (1) -- (5);
				\draw (2) -- (4);
			\end{tikzpicture}%
		}
	\end{minipage}%
	\hspace{1cm}
	\begin{minipage}[b]{0.2\textwidth}
		\centering
		\resizebox{\linewidth}{!}{%
			\begin{tikzpicture}[auto, every node/.style={circle, draw, fill=black,
					inner sep=5pt}]
				\node[fill=black, circle, label={}] (0) at (60:3cm)  {};
				\node[fill=black, circle, label={}] (1) at (120:3cm) {};
				\node[fill=black, circle, label={}] (2) at (180:3cm) {};
				\node[fill=black, circle, label={}] (3) at (240:3cm) {};
				\node[fill=black, circle, label={}] (4) at (300:3cm) {};
				\node[fill=black, circle, label={}] (5) at (0:3cm)   {};

				\draw (0) -- (4);
				\draw (0) -- (5);
				\draw (0) -- (3);
				\draw (1) -- (4);
				\draw (2) -- (5);
				\draw (4) -- (5);
			\end{tikzpicture}%
		}
	\end{minipage}
	\caption{A pair of coflow graphs with different numbers of edges.}
	\label{coflow_graphs_with_diff_edges}
\end{figure}

\begin{remark}
	Coflow mates need not have the same number of edges. The graphs shown in
	Figure~\ref{coflow_graphs_with_diff_edges} form the smallest such pair.
\end{remark}

\begin{table}[htbp]
	\caption{Connected graphs belonging to non-singleton coflow classes on up to 10 vertices.}
	\label{tab2}
	\begin{tabular*}{\tblwidth}{@{}lcccccc@{}}
		\toprule
		& \multicolumn{6}{c}{$n$} \\
		\cmidrule(lr){2-7}
		& \textbf{5} & \textbf{6} & \textbf{7} & \textbf{8} & \textbf{9} & \textbf{10} \\
		\midrule
		\textbf{\# connected graphs} & 21 & 112 & 853 & 11,117 & 261,080 & 11,716,571 \\
		\textbf{\# coflow graphs}    & 0  & 8   & 13  & 63     & 737     & 730,276    \\
		\bottomrule
	\end{tabular*}
\end{table}

\section{Cominimal and Comaximal Graphs}\label{sec5}

Two graphs are said to be \emph{cominimal} if they share the same least eigenvalue of the
average flow matrix. Similarly, two graphs are said to be \emph{comaximal} if they share
the same largest eigenvalue. Since cospectral graphs necessarily have the same minimum and
maximum eigenvalues, the classes of cominimal and comaximal graphs provide natural filters
for the enumeration of cospectral graphs.

\begin{figure}[htbp]
	\centering
	\begin{subfigure}{0.45\linewidth}
		\centering
		\includegraphics[width=\linewidth]{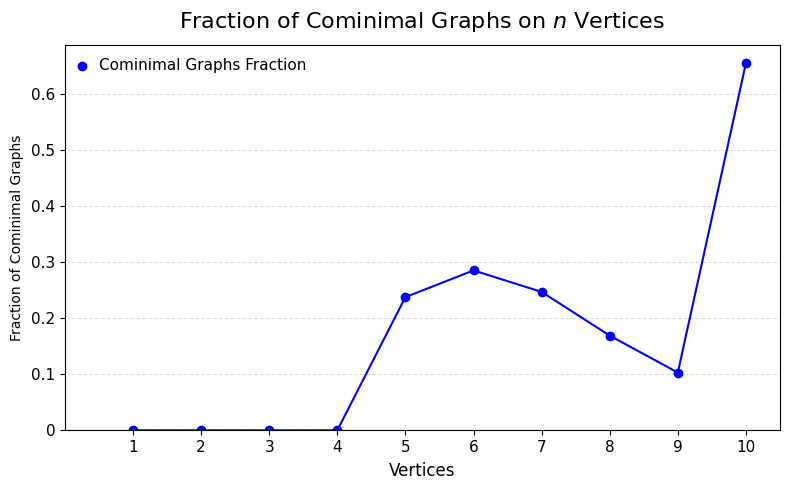}
		\caption{Fraction of cominimal graphs on $n$ vertices.}
		\label{fig:comin_graphs}
	\end{subfigure}
	\hfill
	\begin{subfigure}{0.45\linewidth}
		\centering
		\includegraphics[width=\linewidth]{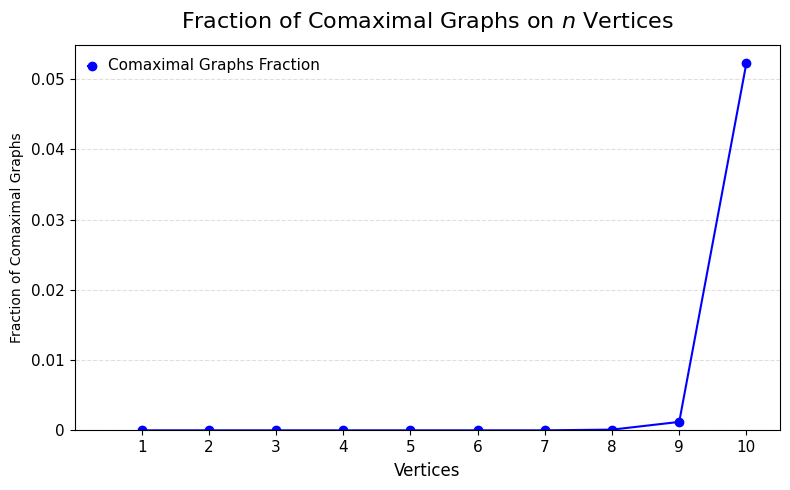}
		\caption{Fraction of comaximal graphs on $n$ vertices.}
		\label{fig:comax_graphs}
	\end{subfigure}
	\caption{Fractions of graphs having a cominimal mate or a comaximal mate on $n$ vertices.}
	\label{fig:comin_max_graphs_frac}
\end{figure}

We therefore enumerate, for each order $n$, the connected graphs having at least one
cominimal mate and the connected graphs having at least one comaximal mate. The resulting
fractions are summarized in Table~\ref{tab:comin-comax-summary}. The proportion of graphs
having a cominimal mate varies substantially with $n$: it is approximately $23\%$ for
$n=5$, increases slightly to $28\%$ for $n=6$, and then decreases to $24.73\%$, $16.92\%$,
and $10.32\%$ for $n=7,8,9$, respectively, as shown in
Figure~\ref{fig:comin_graphs}. On 10 vertices, however, this fraction rises sharply to
approximately $65.62\%$.

\begin{table}[htbp]
	\caption{Fractions of connected graphs having a cominimal mate, a comaximal mate, and
	additional coflow constraints.}
	\label{tab:comin-comax-summary}
	\centering
	\begin{tabular}{ccccc}
		\toprule
		$n$ & Cominimal & Comaximal & Coflow + Cominimal & Coflow + Comaximal \\
		\midrule
		5  & 0.2300 & 0      & --      & 0 \\
		6  & 0.2800 & 0      & 0.0178  & 0 \\
		7  & 0.2473 & 0      & 0.0023  & 0 \\
		8  & 0.1692 & 0.0001 & 0.0025  & 0 \\
		9  & 0.1032 & 0.0012 & 0.0009  & 0 \\
		10 & 0.6562 & 0.0523 & 0.0001  & 0 \\
		\bottomrule
	\end{tabular}
\end{table}

\begin{figure}[htbp]
	\centering
	\begin{subfigure}{0.45\linewidth}
		\centering
		\includegraphics[width=\linewidth]{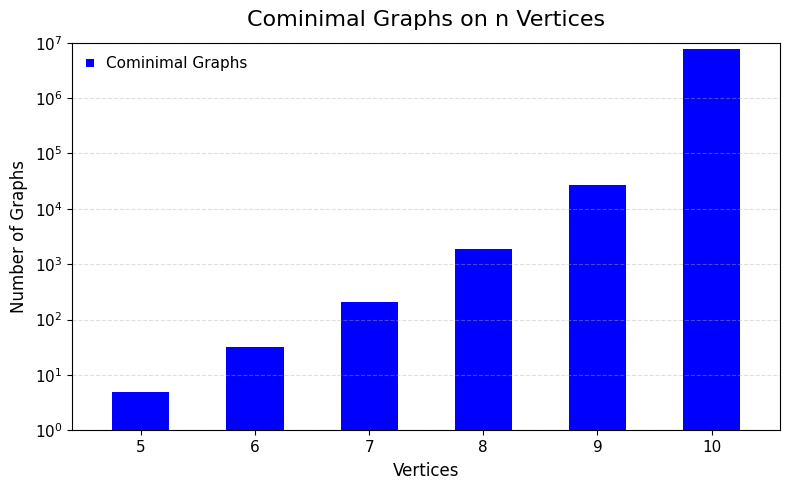}
		\caption{Number of cominimal graphs on $n$ vertices.}
		\label{fig:comin_histogram}
	\end{subfigure}
	\hfill
	\begin{subfigure}{0.45\linewidth}
		\centering
		\includegraphics[width=\linewidth]{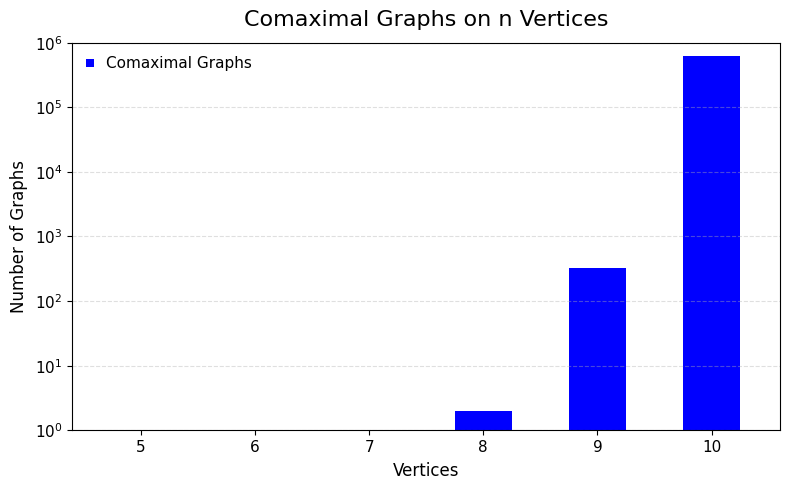}
		\caption{Number of comaximal graphs on $n$ vertices.}
		\label{fig:comax_histogram}
	\end{subfigure}
	\caption{Numbers of cominimal and comaximal graphs on $n$ vertices.}
	\label{fig:comin_max_histograms}
\end{figure}

By contrast, very few graphs have a comaximal mate. Up to 7 vertices, no connected graph
has a comaximal mate. On 8, 9, and 10 vertices, the corresponding fractions are
approximately $0.01\%$, $0.12\%$, and $5.23\%$, respectively, as shown in
Figure~\ref{fig:comax_graphs}. Figure~\ref{fig:comax_histogram} shows that the number of
comaximal graphs increases rapidly at $n=10$.

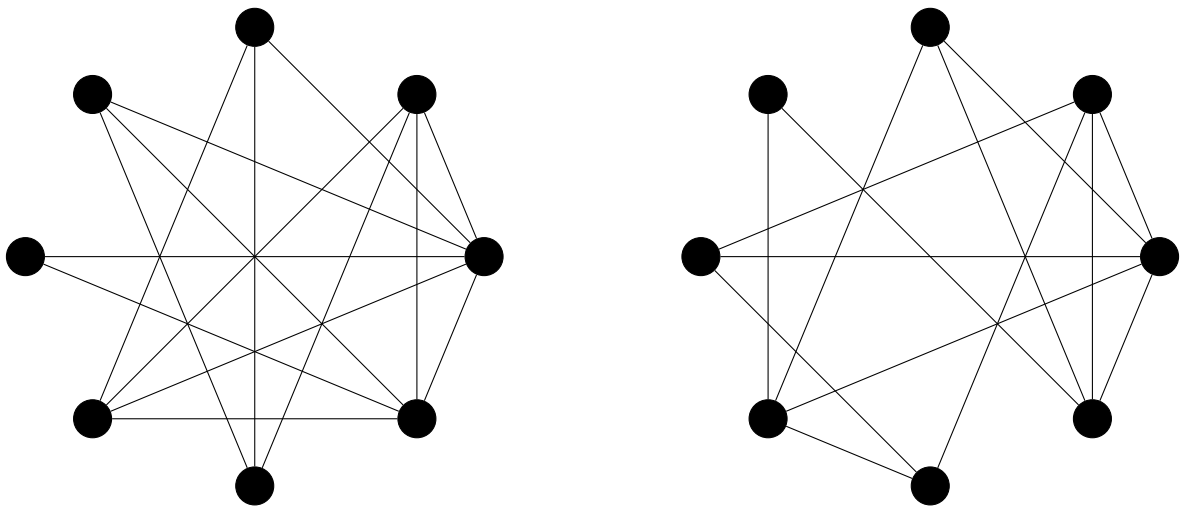
\begin{figure}[htbp]
	\centering
	\resizebox{0.2\textwidth}{!}{%
		\begin{tikzpicture}[auto, every node/.style={circle, draw, fill=black, inner sep=5pt}]
			\node[fill=black, circle, label={}] (0) at (45:3cm)  {};
			\node[fill=black, circle, label={}] (1) at (90:3cm)  {};
			\node[fill=black, circle, label={}] (2) at (135:3cm) {};
			\node[fill=black, circle, label={}] (3) at (180:3cm) {};
			\node[fill=black, circle, label={}] (4) at (225:3cm) {};
			\node[fill=black, circle, label={}] (5) at (270:3cm) {};
			\node[fill=black, circle, label={}] (6) at (315:3cm) {};
			\node[fill=black, circle, label={}] (7) at (360:3cm) {};

			\draw (0) -- (4); \draw (0) -- (5);
			\draw (0) -- (6); \draw (0) -- (7);
			\draw (1) -- (4); \draw (1) -- (5); \draw (1) -- (7);
			\draw (2) -- (5); \draw (2) -- (6); \draw (2) -- (7);
			\draw (3) -- (6); \draw (3) -- (7);
			\draw (4) -- (6); \draw (4) -- (7);
			\draw (6) -- (7);
		\end{tikzpicture}%
	}
	\hspace{1cm}
	\resizebox{0.2\textwidth}{!}{%
		\begin{tikzpicture}[auto, every node/.style={circle, draw, fill=black, inner sep=5pt}]
			\node[fill=black, circle, label={}] (0) at (45:3cm)  {};
			\node[fill=black, circle, label={}] (1) at (90:3cm)  {};
			\node[fill=black, circle, label={}] (2) at (135:3cm) {};
			\node[fill=black, circle, label={}] (3) at (180:3cm) {};
			\node[fill=black, circle, label={}] (4) at (225:3cm) {};
			\node[fill=black, circle, label={}] (5) at (270:3cm) {};
			\node[fill=black, circle, label={}] (6) at (315:3cm) {};
			\node[fill=black, circle, label={}] (7) at (360:3cm) {};

			\draw (0) -- (3); \draw (0) -- (5);
			\draw (0) -- (6); \draw (0) -- (7);
			\draw (1) -- (4); \draw (1) -- (6); \draw (1) -- (7);
			\draw (2) -- (4); \draw (2) -- (6);
			\draw (3) -- (5); \draw (3) -- (7);
			\draw (4) -- (5); \draw (4) -- (7);
			\draw (6) -- (7);
		\end{tikzpicture}%
	}
	\caption{The unique pair of comaximal graphs on 8 vertices.}
	\label{fig:comaximal_8vertices}
\end{figure}

Although cominimality and comaximality provide valid spectral filters, they are not by
themselves sufficiently restrictive for the final cospectral search. For instance, on
10 vertices there are $7{,}688{,}945$ connected graphs having a cominimal mate and
$613{,}212$ having a comaximal mate, leaving a search space that is still computationally
large.

\begin{figure}[htbp]
	\centering
	\begin{subfigure}{0.45\linewidth}
		\centering
		\includegraphics[width=\linewidth]{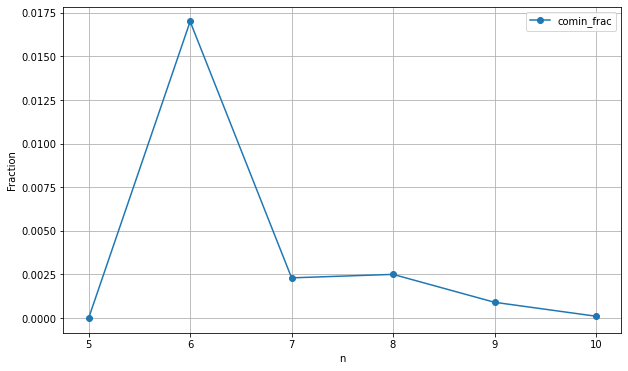}
		\caption{Fraction of cominimal graphs among coflow graphs on $n$ vertices.}
		\label{fig:coflow_comin_graphs}
	\end{subfigure}
	\hfill
	\begin{subfigure}{0.45\linewidth}
		\centering
		\includegraphics[width=\linewidth]{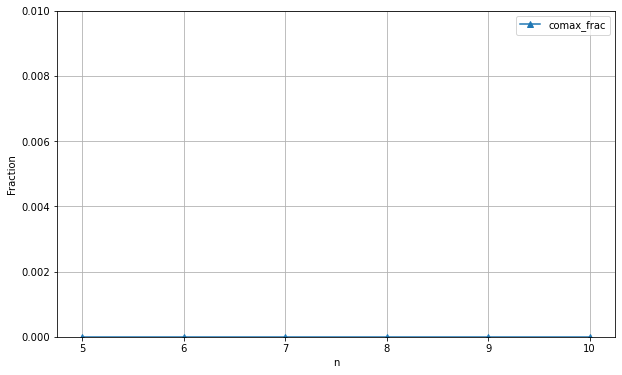}
		\caption{Fraction of comaximal graphs among coflow graphs on $n$ vertices.}
		\label{fig:coflow_comax_graphs}
	\end{subfigure}
	\caption{Fractions of coflow graphs that are also cominimal or comaximal on $n$ vertices.}
	\label{fig:coflow_comin_max_graphs_frac}
\end{figure}

We therefore combine these filters with the coflow condition introduced in
Section~\ref{subsec:coflow}. Since cospectral graphs must be coflow, cominimal, and
comaximal, it is natural to examine the intersections of these classes. The results show
that the fraction of \emph{coflow graphs} that also have a cominimal mate decreases
substantially for larger $n$: \textcolor{black}{ it is $1.78\%$ on $6$ vertices, $0.23\%$ on $7$ vertices,  approximately $0.25\%$ for $n=8$, $0.09\%$ for $n=9$,}
and only $0.01\%$ for $n=10$, as shown in Figure~\ref{fig:coflow_comin_graphs}. Moreover,
up to 10 vertices, no coflow graph has a comaximal mate, as shown in
Figure~\ref{fig:coflow_comax_graphs}.

These observations indicate that the coflow condition is considerably more effective than
the extremal-eigenvalue filters in reducing the search space. Consequently, it is
computationally preferable to first restrict attention to coflow classes and then, if
necessary, refine the search further using comaximality.

\begin{figure}[htbp]
	\centering
	\begin{subfigure}{0.45\linewidth}
		\centering
		\includegraphics[width=\linewidth]{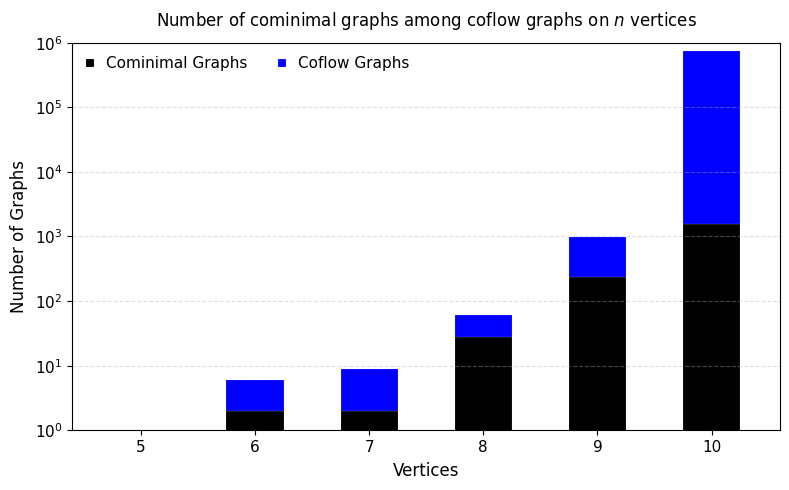}
		\caption{Number of cominimal graphs among coflow graphs on $n$ vertices.}
		\label{fig:coflow-comin_histogram}
	\end{subfigure}
	\hfill
	\begin{subfigure}{0.45\linewidth}
		\centering
		\includegraphics[width=\linewidth]{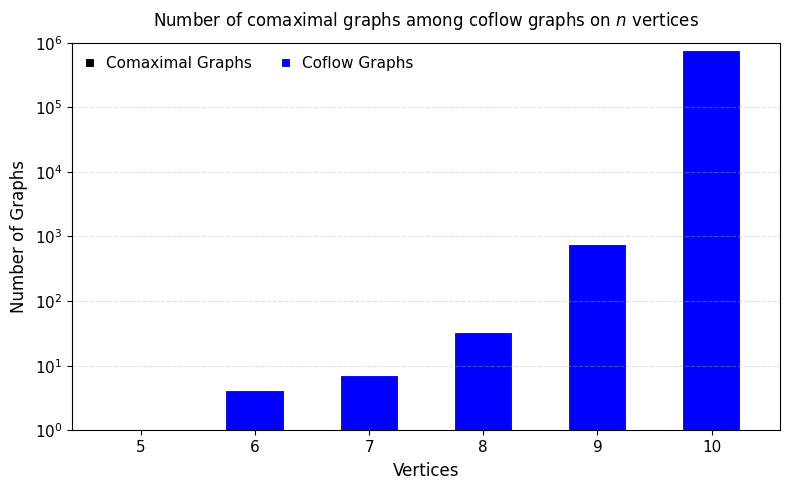}
		\caption{Number of comaximal graphs among coflow graphs on $n$ vertices.}
		\label{fig:coflow-comax_histogram}
	\end{subfigure}
	\caption{Cominimal and comaximal graphs among coflow graphs on $n$ vertices.}
	\label{fig:comin_comax_coflow_histograms}
\end{figure}

\section{Enumeration Results}\label{sec6}

We now present the main enumeration results for the average flow matrix and compare them
with the best previously reported results for other graph-associated matrices. Godsil and
McKay~\cite{r} enumerated cospectral graphs with respect to the adjacency matrix $A$ for
graphs on up to 9 vertices. Haemers and Spence~\cite{f} extended this line of work to the
matrices $A$, $L$, and $Q$, obtaining the strongest results for $Q$. Wilson and
Zhu~\cite{y} showed that improved enumeration can be achieved by considering several
matrices simultaneously, namely $A$, $L$, $Q$, and $\tilde{L}$. Aouchiche and
Hansen~\cite{d} further improved the enumeration results by considering the
distance-associated matrices $D$, $D^L$, and $D^Q$. More recently, Abiad and
Alfaro~\cite{b} studied cospectral and coinvariant graphs, and showed that additional
improvement can be obtained by combining adjacency- and distance-based invariants.

For each $n\leq 10$, we generated all connected graphs on $n$ vertices and computed their
average flow matrices. 
We restricted the enumeration to connected graphs on at most 10 vertices because the
number of connected graphs grows extremely rapidly with $n$, making exhaustive spectral
comparisons beyond this range computationally prohibitive on standard hardware.
To reduce the computational cost of exhaustive pairwise spectral
comparisons, we first partitioned the graphs into coflow classes, since two graphs can be
cospectral with respect to $F$ only if they have the same average flow sum. We then used
the largest eigenvalue as an additional filter, retaining only those graphs that were both
coflow and comaximal. Finally, for the remaining candidates, we compared the full spectra
of their average flow matrices numerically in order to identify cospectral mates.

All computations were performed on a machine equipped with an
Intel\textsuperscript{\textregistered} Core\textsuperscript{\texttrademark} i5-8365U CPU @
1.60\,GHz, 8.00\,GB RAM, a 64-bit operating system, and Python~3.11. The eigenvalues were
computed numerically using 7-digit precision, so that the resulting counts are directly
comparable with those reported in the existing literature. The results are summarized in
Table~\ref{enum3}.

\begin{table}[htbp]
	\caption{Number of connected graphs that have a cospectral or cospectral-coinvariant mate
	with respect to various graph-associated matrices.}
	\label{enum3}
	\begin{tabular*}{\tblwidth}{@{}p{5cm}p{3.8cm}cccccc@{}}
		\toprule
		\textbf{Reference} & \textbf{Family} & \multicolumn{6}{c}{$n$} \\
		\cmidrule(lr){3-8}
		& & \textbf{5} & \textbf{6} & \textbf{7} & \textbf{8} & \textbf{9} & \textbf{10} \\
		\midrule
		\cite{r}
		& $\mathcal{G}_{n}^{cs}(A)$
		& 0 & 2 & 63 & 1,353 & 46,930 & NA \\
		\cite{f}
		& $\mathcal{G}_{n}^{cs}(Q)$
		& 2 & 10 & 80 & 1,047 & 17,627 & 615,919 \\
		\cite{y}
		& $\mathcal{G}_{n}^{cs}(A, L, Q, \tilde{L})$
		& 0 & 0 & 0 & 0 & 2 & 10,124 \\[4pt]
		\multirow{2}{*}{\cite{d}}
		& $\mathcal{G}_{n}^{cs}(D, D^Q)$
		& 0 & 0 & 0 & 0 & 0 & 7,712 \\
		& $\mathcal{G}_{n}^{cs}(D, D^L, D^Q)$
		& 0 & 0 & 0 & 0 & 0 & 7,622 \\[4pt]
		\multirow{2}{*}{\cite{b}}
		& $\mathcal{G}_{n}^{cs-in-in}(D, D^L, D^Q)$
		& 0 & 0 & 0 & 0 & 0 & 7,510 \\
		& $\mathcal{G}_{n}^{cs-in-in}(A, D^L, D^Q)$
		& 0 & 0 & 0 & 0 & 0 & 7,490 \\[4pt]
		Our results
		& $\mathcal{G}_{n}^{cs}(F)$
		& \textbf{0} & \textbf{0} & \textbf{0} & \textbf{0} & \textbf{0} & \textbf{0} \\
		\bottomrule
	\end{tabular*}
\end{table}
\FloatBarrier

Table~\ref{enum3} shows that, among the graph-associated matrices considered so far in the
literature, the average flow matrix provides the strongest spectral discrimination on
connected graphs with up to 10 vertices. In particular, within the adopted numerical
precision, our computations yielded no connected graph on at most 10 vertices having a
cospectral mate with respect to $F$. This improves upon all previously reported results
listed in Table~\ref{enum3}, for which nontrivial cospectral classes still remain at
order~10.

		\section{Conclusions}\label{sec7}

In this paper, we introduced the {average flow matrix}, a new graph-associated matrix
that combines the average lengths of all simple paths originating at each vertex with the
pairwise shortest-path distances between vertices. The primary aim of this work was to
assess the effectiveness of this matrix in the enumeration of cospectral graphs.

In addition, we established several supporting properties of the average flow matrix,
including its behavior on flow-stable graphs and complete graphs, and introduced the
notions of coflow, cominimal, and comaximal graphs as natural spectral filters. While
these results are primarily supportive in nature, they help clarify the structure of the
matrix and provide useful tools for the computational analysis.

The main contribution of the paper is computational. We carried out an exhaustive
enumeration of connected graphs on up to 10 vertices and, within the adopted numerical
precision, found no connected graph having a cospectral mate with respect to the average
flow matrix. In this sense, the average flow matrix provides stronger spectral
discrimination on this range than the graph-associated matrices considered in the existing
literature, including those studied by Godsil and McKay~\cite{r}, Haemers and
Spence~\cite{f}, Wilson and Zhu~\cite{y}, Aouchiche and Hansen~\cite{d}, and Abiad and
Alfaro~\cite{b}. In particular, it improves on the best previously reported enumeration
results, all of which still leave nontrivial cospectral classes at order 10.

Our experiments also yield several useful observations about the search space. Among the
filters considered in this paper, the coflow condition provides the strongest reduction,
while comaximality offers a useful secondary refinement and cominimality is comparatively
less effective. These observations suggest that the average flow sum captures a substantial
amount of spectral information and is particularly valuable in large-scale enumeration.

The restriction to graphs on at most 10 vertices reflects the rapid growth in the number of
connected graphs and the resulting computational cost of exhaustive spectral comparison.
Several directions for future work remain open. 
It would be natural to extend the enumeration beyond 10 vertices using improved filtering
strategies or greater computational resources, and to investigate exact verification
methods for cospectrality with respect to the average flow matrix. Another interesting
direction is to study how the matrix behaves on important graph classes such as trees,
distance-regular graphs, and other highly symmetric families.
 A more detailed
analysis of coflow, cominimal, and comaximal classes may also help explain why the average
flow matrix performs so well in practice.

		\bibliographystyle{cas-model2-names}
		
		\bibliography{cas-refs}
		
		
		
	\end{document}